\documentclass{elsart}
 \journal{Applied Mathematics and Computation}

\newtheorem{proposition}{Proposition}[section]

\newtheorem{remark}{Remark}[section]

\newtheorem{example}{Example}[section]

\def \IR {{\rm I \! R }}

\def \N {{\rm I  \! N }}

\begin{document}

\begin{frontmatter}
\title{ Classes of  second order
 nonlinear differential equations reducible to first order ones by variation of parameters}
\author[cotonou]{Mahouton Norbert Hounkonnou\thanksref{XX}}and
\author[cotonou]{Pascal Alain Dkengne Sielenou}

\thanks[XX]{Correspondence should be addressed to:
 norbert.hounkonnou@cipma.uac.bj, with copy to hounkonnou@yahoo.fr}
\address[cotonou]{International Chair in Mathematical Physics
and Applications (ICMPA-UNESCO Chair),\\ 072 B.P.:50,  Cotonou,
Republic of Benin}

\begin{abstract}
The method of paramùeter variation for linear differential equations
is extended to classes of second order nonlinear differential
equations. This allows to reduce the latter to first order
differential equations. Known classical equations such as the
Bernoulli, Riccati and Abel equations are recovered in illustrated
relevant examples.
\end{abstract}
\end{frontmatter}

\section{Introduction}
Second order nonlinear differential equations are important
 for investigation of nonlinear phenomena in different fields
 of physics and mathematics. They are used to model a wide number of phenomena in
 plasma physics, solid state physics, optics, bio-hydrodynamics,
 chemical processes, nonlinear quantum mechanics, etc. See [\ref{r23}-\ref{hk}]
 and references therein. Formal methods
of analytical integration of nonlinear differential equations thus
appear  of great interest in the theory of differential equations.
To cite a few, one can mention various powerful methods, which are
most familiarly  used, such as group symmetry methods
[\ref{r20}-\ref{r21}], tanh method [\ref{r9}-\ref{r11}], extended
tanh method [\ref{r12}-\ref{r14}], sine-cosine method
[\ref{r15}-\ref{r16}], Jacobi elliptic method [\ref{r17}-\ref{r18}],
Backlund transformations [\ref{r25}-\ref{r26}], inverse scattering
method [\ref{r27}], pseudo spectral method [\ref{r28}] and
F-expansion method [\ref{r29}-\ref{r31}].

 Unfortunately, only a few types of nonlinear
differential equations can be exactly solved, what explains the
permanent need to recur to novel
 tools of handling differential equations. In such a direction,
  the method of reducing the order of ordinary differential equations (ODEs)
can be of inestimable usefulness in most cases.
 Indeed,
  in general, the lowest the order of a differential equation is,
 the greatest is the possibility of finding  suitable analytical methods for its resolution.
 Moreover, the accuracy of  numerical
treatment of a differential equation decreases as the order of the
equation increases. Therefore, when a differential equation cannot
be directly analytically solved, it is meaningful to transform it
into an equation of lower order. The probability of solving the
 reduced equation is greater than the probability of
solving  the original equation. In practice, it is generally
difficult to perform such a reduction.

 This paper aims at proving that the  usual method of
parameter variation  can be successfully used in some cases to
reduce by one the order of  nonlinear  differential equations
(NLDEs). In the sequel, we provide with  some classes of second
order nonlinear differential equations,  transformable to first
order ones using specific parameter variations. Concrete examples
are exhibited as matter of illustration. Specifically, the classes
of equations of the forms
\begin{equation}\label{eq1}
\left(y'\right)^{m}y''+a(x)\left(y'\right)^{m+1}=f(x,y,y'),
\end{equation}
and
\begin{equation}\label{eq2}
\left(y' \right)^m y''+a(y')(y')^{m+2}=f(x,y,y'),
\end{equation}
where $m$ is a positive integer and $a$ is an integrable function,
are considered. Note that, if the function $f(x,y,y')=g(x,y'),$ i.e.
the second members of equations (\ref{eq1}) and (\ref{eq2}) do not
explicitly depend  on the variable $y$, then we can perform  the
natural change of variables $z(x)=y'(x)$ to lower the order of these
equations. Besides, if the function $a$   is constant only in
equation (\ref{eq1}), (not necessarily constant  in equation
(\ref{eq2})), and if $f(x,y,y')=g(y,y'),$ i.e. the equations
(\ref{eq1}) and (\ref{eq2}) are autonomous, then we can perform the
 change of variables $w(y)=y'(x)$ to transform  the above mentioned
 classes of
equations into first order ones.


\section{First class of reducible second order  NLDEs}
Let us consider nonlinear second order differential equations  of
the following type:
\begin{equation}\label{eq3}
y''+a(x)y'=F(x)+y'G(y)e^{-\int a(x) dx}
\end{equation}
whose the linear  part, i.e.
\begin{equation}\label{eq4}
y''+a(x)y'=0
\end{equation}
yields the solution
\begin{equation}\label{eq5}
y'=Ce^{-\int a(x) dx},
\end{equation}
where $C$ is an arbitrary constant. Suppose that $C$ is a
differentiable function of both variables $x$ and $y$ expressed in
the form
\begin{equation}\label{eq6}
C=H(x)+K(y).
\end{equation}
Then, (\ref{eq5}) can be rewritten as
\begin{equation}\label{eq7}
y'=\left[H(x)+K(y)\right]e^{-\int a(x) dx}
\end{equation}
that we differentiate  to obtain
\begin{equation}\label{eq8}
y''=e^{-\int a(x) dx}\left[H'_{x}+y'K'_y\right] -a(x)y',
\end{equation}
where
 $$H'_x=\frac{dH(x)}{dx} \quad\mbox{and}\quad
K'_y=\frac{dK(y)}{dy}.$$
Substituting (\ref{eq7}) and (\ref{eq8}) into (\ref{eq3}), we find
\begin{equation}\label{eq9}
H'_{x}+y'K'_y= F(x)e^{\int a(x) dx} +y'G(y).
\end{equation}
Clearly, (\ref{eq9}) will take place if
\begin{equation}\label{eq10}
H'_{x}=F(x)e^{\int a(x) dx}\quad \mbox{and}\quad K'_y=G(y).
\end{equation}
We therefore state the following result:
\begin{proposition}\label{pn1}
The  second order nonlinear differential equation (\ref{eq3}) can be
reduced to the first order differential equation (\ref{eq7}), where
the functions $H$ and $K$ are solutions of the first order
differential equations (\ref{eq10}),  respectively.
\end{proposition}
As a matter of clarity, let us consider the following example:
\begin{example}
\emph{Let us consider the function $G$ in (\ref{eq3}) in the form}
\begin{equation}
 G(y)=\frac{d}{dy}\left(
\frac{b_0+b_1y+b_2y^2+b_3y^3}{c_0+c_1y} \right),
\end{equation}
\emph{where} $b_0,\,b_1,\,b_2,\,b_3$ \emph{are arbitrary constants
and} $c_0,\,c_1$ \emph{are constants such that} $(c_0,c_1) \ne
(0,0).$ \emph{Then equations (\ref{eq10}) yield :}
\begin{eqnarray}
H(x)&=&\int F(x) e^{\int a(x)dx}dx \cr
K(y)&=&A+\frac{b_0+b_1y+b_2y^2+b_3y^3}{c_0+c_1y},
\end{eqnarray}
\emph{where} $A$ \emph{is an arbitrary constant of integration.}
\emph{Therefore, equation (\ref{eq7}) is reduced to the well known
Abel equation of second kind}
\begin{eqnarray}\label{eq11}
y'=e^{-\int
a(x)dx}\,\frac{b_0+c_0(H(x)+A)+\left[c_1(H(x)+A)+b_1\right]y+b_2y^2+b_3y^3}{c_0+c_1y}.\nonumber\\
\end{eqnarray}
\emph{Finally, the equation (\ref{eq3}) takes the form}
\begin{equation}
y''+a(x)y'=F(x)+y'\frac{d}{dy}\left(
\frac{b_0+b_1y+b_2y^2+b_3y^3}{c_0+c_1y} \right)e^{-\int a(x) dx}
\end{equation}
\emph{ and is integrable if the constants} $b_0,\,b_1,\,b_2,\,b_3,$
$c_0,\,c_1$ \emph{and the functions} $a$ \emph{and} $F$ \emph{are
chosen in such a way  that the Abel equation (\ref{eq11}) be
integrable.}

\emph{In particular, for} $c_0=1,\,c_1=0$ \emph{and} $F=0,$
\emph{equation (\ref{eq11}) leads to the separable Abel equation}
\begin{equation}
y'=(A+b_0+b_1y+b_2y^2+b_3y^3)\,e^{-\int a(x) dx}
\end{equation}
\emph{whose the implicit solution is given by}
\begin{equation}
\int\frac{dy}{A+b_0+b_1y+b_2y^2+b_3y^3}=\int e^{-\int a(x) dx}dx +B,
\end{equation}
\emph{where} $B$ \emph{is an arbitrary constant of integration.}
\end{example}

It is worth noticing that   a second
kind Abel equation of the form
\begin{equation}\label{eq12}
y'=\frac{f_3y^3+f_2y^2+f_1y+f_0}{g_1y+g_0},\quad \mbox{with}\quad
f_3\neq 0,
\end{equation}
where $f_i, \,(i=0,\,1,\,2,\,3),$ and $g_j,\,(j=0,\,1),$ are
arbitrary functions of $x$, can be transformed into a canonical
form. Indeed, using the variable change
\begin{equation}
\left\{ x=t,\,\,y=\frac{1-g_0u}{g_1u}  \right\},
\end{equation}
where $t$ and $u=u(t)$ are  the new independent and
dependent variables, respectively, equation (\ref{eq12}) becomes
\begin{equation}\label{eq13}
u'_t=\widetilde{f}_3u^3+\widetilde{f}_2u^2+\widetilde{f}_1u+\widetilde{f}_0.
\end{equation}
Making use of the substitution
\begin{equation}
u=v-\frac{\widetilde{f}_2}{3\,\widetilde{f}_3},
\end{equation}
equation (\ref{eq13}) can be put in the form
\begin{equation}
v'_t=h_3v^3+h_1v+h_0.
\end{equation}
Now, setting
\begin{equation}
v=E(t)w,\quad \mbox{where}\quad E(t)=e^{\int h_1(t) dt},
\end{equation}
brings this equation to the simpler form:
\begin{equation}\label{eq14}
w'_t=\widetilde{h}_3w^3+\widetilde{h}_0,
\end{equation}
which, in turn,  can be reduced, with the help of the new
independent variable
\begin{equation}
s=\int \widetilde{h}_3(t) dt,
\end{equation}
to the usual  canonical form of  Abel equation of the first kind
\begin{equation}\label{eq15}
w'_s=w^3(s)+k(s).
\end{equation}
The latter is integrable by various methods known in the literature.
See [\ref{r2}] and [\ref{r5}] (and references therein) for a good
compilation of techniques developed to solve (\ref{eq15}) for
particular expressions  of $k(s)$.


\section{Second class of reducible second order  NLDEs}
In this section, we discuss the second relevant type of nonlinear second order differential
 equations compilable
in the following general form:
\begin{equation}\label{eqb3}
y''+a(y)(y')^{2}=F(x)e^{-\int a(y) dy}+ y'G(y).
\end{equation}
We first solve the left hand side  part of this equation,
namely
\begin{equation}\label{eqb4}
y''+a(y)(y')^{2}=0
\end{equation}
to get
\begin{equation}\label{eqb5}
y'=Ce^{-\int a(y) dy},
\end{equation}
where $C$ is an arbitrary constant. As in the previous section, we
assume that $C$ is a differentiable function of both $x$ and $y$ as
follows
\begin{equation}\label{eqb6}
C=H(x)+K(y).
\end{equation}
Then, the equation (\ref{eqb5}) takes the form
\begin{equation}\label{eqb7}
y'=\left[H(x)+ K(y)\right]e^{-\int a(y) dy}
\end{equation}
that we differentiate to obtain
\begin{equation}\label{eqb8}
y''=e^{-\int a(y) dy}\left[H'_{x}+y'K'_y\right] -a(y)(y')^{2}.
\end{equation}
Substituting (\ref{eqb7}) and (\ref{eqb8}) into (\ref{eqb3}), we
find
\begin{equation}\label{eqb9}
H'_{x}+y'K'_y= F(x)+ y'G(y)e^{\int a(y) dy}.
\end{equation}
Clearly, (\ref{eqb9}) will take place if
\begin{equation}\label{eqb10}
H'_{x}=F(x)\quad \mbox{and}\quad K'_y=G(y)e^{\int a(y) dy}.
\end{equation}
Therefore, the following statement holds:
\begin{proposition}\label{pn2}
The  second order nonlinear differential equation (\ref{eqb3}) can
be reduced to the first order differential equation (\ref{eqb7}),
where the functions $H$ and $K$ are   solutions of the two first
order differential equations (\ref{eqb10}), respectively.
\end{proposition}
For illustration, let us consider the  following example.
\begin{example}
\emph{Let} $a(y)=-\frac{1}{y}$ \emph{and} $G(y)=\beta\,y^n,$
\emph{where} $n$ \emph{is a non zero positive integer and} $\beta$
\emph{is a constant. Then equation (\ref{eqb3}) becomes}
\begin{equation}\label{eqb11}
y''-\frac{1}{y}(y')^{2}=F(x)y+ \beta\,y'y^n.
\end{equation}

\emph{Equations (\ref{eqb10}) yield}
\begin{equation}
H(x)=\int F(x)dx+A\quad \mbox{\emph{and}}\quad
K(y)=\frac{\beta}{n}\,y^n,
\end{equation}
\emph{where} $A$ \emph{is an arbitrary constant of integration.}

\emph{Therefore, equation (\ref{eqb7}) becomes the Bernoulli
equation [\ref{r7}]}
\begin{equation}\label{bern}
y'=(H(x)+A)\,y+\frac{\beta}{n}\,y^{n+1}.
\end{equation}

\emph{The substitution} $w(x)=y^{1-n}$ \emph{transforms (\ref{bern})
into the linear equation}
\begin{equation}
w'_x=-n\,(H(x)+A)\,w-\beta
\end{equation}
\emph{whose the solution is}
\begin{equation}
w(x)=\frac{-\beta \int e^{n\int (H(x)+A)\,dx} dx+B}{e^{n\int
(H(x)+A)\,dx}},
\end{equation}
\emph{where} $B$ \emph{is an arbitrary constant of integration.}
\end{example}


\section{Third class of reducible second order  NLDEs}
The third group of  second order nonlinear differential equations
 can be expressed as
\begin{equation}\label{eqc1}
\left(y' \right)^{m}y''+a(x)\left(y'\right)^{m+1}=e^{-(m+2)\int
a(x)dx}F\left( y,\,y'e^{\int a(x)dx}\right).
\end{equation}
\begin{remark}
\emph{If we set in (\ref{eq3})} $F(x)=C\,e^{-2\int a(x)dx},$
\emph{where} $C$ \emph{is a constant, then equation (\ref{eq3})
appears as a particular case of equations (\ref{eqc1}) considered
with} $m=0.$
\end{remark}
The linear  part of equation (\ref{eqc1})
\begin{equation}\label{eqc2}
\left(y' \right)^{m}y''+a(x)\left(y'\right)^{m+1}=0
\end{equation}
can be readily solved to give
\begin{equation}\label{eqc3}
\left(y' \right)^{m+1}=K^{m+1}e^{-(m+1)\int a(x) dx},
\end{equation}
where $K$ is an arbitrary constant. Now, suppose that $K$ is a
differentiable function of the variable $y.$ Then, the equation
(\ref{eqc3}) becomes
\begin{equation}\label{eqc4}
\left(y' \right)^{m+1}=K^{m+1}(y)e^{-(m+1)\int a(x) dx}.
\end{equation}
Differentiate (\ref{eqc4}) to obtain
\begin{equation}\label{eqc5}
\left(y' \right)^{m}y^{(n+2)}=K^{m+1}(y)K'_{y}e^{-(m+2)\int a(x) dx}
-a(x)\left(y' \right)^{m+1}.
\end{equation}
Substituting (\ref{eqc4}) and (\ref{eqc5}) into (\ref{eqc1}), we
find
\begin{equation}\label{eqc6}
K^{m+1}(y)K'_{y}=F\left( y,\,K(y)\right).
\end{equation}
We therefore set the following result:
\begin{proposition}\label{pn3}
The  second order nonlinear differential equation (\ref{eqc1}) can
be reduced to the first order differential equation (\ref{eqc6})
using (\ref{eqc4}). Furthermore, if $K(y)=\Phi(y,\,A),$ where $A$ is
an arbitrary constant, is the general solution of (\ref{eqc6}), then
the general solution of (\ref{eqc1}) is given by
\begin{equation}\label{eqc7e}
\int\frac{dy}{\Phi(y,\,A)}=\int e^{-\int a(x) dx}dx+B,
\end{equation}
where $B$ is an arbitrary constant.
\end{proposition}
The following examples are of particular significance.
\begin{example}
\emph{Suppose in (\ref{eqc1}) that} $F$ \emph{has one of the
following forms}
\begin{eqnarray}
&&i) \,F(u,\,v)=u^{m+1}\,f\left(  \frac{au+bv+c}{\alpha u+\beta
v+\gamma}\right)\cr && \mbox{\emph{or}}\cr && ii)
F(u,\,v)=v^{m+1}\,f\left( \frac{au+bv+c}{\alpha u+\beta
v+\gamma}\right),
\end{eqnarray}
\emph{where} $\alpha,\,\beta,\,\gamma,\,a,\,b$ \emph{and} $c$
\emph{are constants such that the function} $f$ \emph{is well
defined. Then, using the method developed in [\ref{r2}] and
[\ref{r4}], the corresponding first order differential equation
(\ref{eqc6}) can be transformed into   a homogeneous equation which,
in turn, can be reduced to  a separable equation. As a matter of
illustration, let us perform such a transformation in the case
where} $F$ \emph{is of the
 form $ii)$.}
\begin{enumerate}
  \item \emph{If} $\Delta=a\beta-b\alpha\neq 0,$ \emph{the transformations}
  \begin{equation}
  y=u+\frac{b\gamma-c\beta}{\Delta}\quad \mbox{and}\quad K(y)=v(u)+\frac{c\alpha-a\gamma}{\Delta}
  \end{equation}
  \emph{lead to the equation}
  \begin{equation}
  v'_u=f\left( \frac{au+bv}{\alpha u+\beta v} \right).
  \end{equation}
  \emph{Dividing both the numerator and the denominator of the fraction on
  the right-hand side by} $u,$ \emph{we obtain the homogeneous equation}
  \begin{equation}
   v'_u=f\left(\frac{a+b\frac{v}{u}}{\alpha+\beta\frac{v}{u}}\right)\equiv\widetilde{f}\left(
   \frac{v}{u}\right),
  \end{equation}
    \emph{for which the substitution} $w(u)=\frac{v}{u}$ \emph{gives the
    separable equation}
    \begin{equation}
    u\,w'_u=\widetilde{f}(w)-w.
    \end{equation}

  \item \emph{For} $\Delta=0$ \emph{and} $b\neq 0,$ \emph{the substitution}
  $v(y)=ay+bK(y)+c$ \emph{engenders the separable equation}
  \begin{equation}
   v'_y=a+bf\left( \frac{bv}{\beta v+b\gamma-c\beta} \right).
  \end{equation}

  \item \emph{For} $\Delta=0$ \emph{and} $\beta\neq 0,$ \emph{from the substitution} $v(y)=\alpha y+\beta
  K(y)+\gamma,$ \emph{we deduce the separable equation}
  \begin{equation}
   v'_y=\alpha+\beta f\left( \frac{bv+c\beta-b\gamma}{\beta v}
   \right).
  \end{equation}
\end{enumerate}
 \emph{  Setting} $m=0,$ $a(x)=\frac{1}{x}$ \emph{and}
$F(u,\,v)=v\,\left[ \left(  \frac{v}{u}
\right)^2+2\,\left(\frac{v}{u} \right) \right],$ \emph{the equation (\ref{eqc1}) becomes}
\begin{equation}\label{eqx10}
y''+\frac{1}{x}\,y'=x\,y'\,\left( \frac{y'}{y}
\right)^2+2\,\frac{y'^2}{y}.
\end{equation}
\emph{Perform in (\ref{eqx10}) the substitution (\ref{eqc4}) which
takes the form}
\begin{equation}
y'=\frac{K(y)}{x}.
\end{equation}
\emph{By (\ref{eqc6}), the function} $K$  \emph{satisfies the
homogeneous equation}
\begin{equation}\label{eqx11}
K'_y=\left( \left(\frac{K}{y}  \right)^2 +2\,\frac{K}{y} \right)
\end{equation}
\emph{which gives the solution}
\begin{equation}
K(y)=\frac{y^2}{A-y},
\end{equation}
\emph{where} $A$ \emph{is an arbitrary constant of integration.}
\emph{Hence, using the integral (\ref{eqc7e}), we obtain, for
equation (\ref{eqx10}), the implicit  solution}
\begin{equation}
-\frac{A}{y}-\ln (y)=\ln (x)+B ,
\end{equation}
\emph{where} $B$ \emph{is an integration constant.}
\end{example}

Suppose  that $F$, defined in (\ref{eqc1}), has the following general form
\begin{equation}\label{eqc7}
F(u,\,v)=\frac{\sum_{\nu=1}^p\,h_\nu(u)\,F_\nu(v)}{\sum_{\eta=1}^q\,k_\eta(u)\,G_\eta(v)},
\end{equation}
where $p,\,q$ are positive integers, $h_\nu,\,k_\eta,\,F_\nu$ and
$G_\eta$ are arbitrary functions such that $F$ be well defined.
\begin{example}
 \emph{If} $F$ \emph{in (\ref{eqc7}) is of the form}
  \begin{equation}\label{eqx21}
  F(u,\,v)=h_1(u)\,R\left(v,\,\sqrt{P(v)}\right),
  \end{equation}
   \emph{where} $R$ \emph{is a rational function of two variables}  \emph{and the function} $P$ \emph{under the
   radical is a polynomial of  degree  three or four, then the
   corresponding reduced equation (\ref{eqc6}) readily turns to be  an elliptic
   integral which can  be merely integrated by the method presented in
   [\ref{r8}].}
\end{example}
\begin{example}
   \emph{If} $F$ \emph{in (\ref{eqc7}) has  the
  form}
  \begin{equation}
   F(u,\,v)=v^{m+1}(h_1(u)\,v+h_2(u)\,v^n),\quad n\in
   \N\setminus\{0,1\},
  \end{equation}
\emph{then the corresponding reduced equation (\ref{eqc6}) leads to
a Bernoulli equation.}
\emph{  With} $m=0,$
$a(x)={2}/{x}$ \emph{and} $F(u,\,v)=v\,\left( v+v^3 \right),$
\emph{the equation (\ref{eqc1}) becomes}
\begin{equation}\label{eqxx10}
y''+\frac{2}{x}\,y'=(y')^2+\left( x\,y' \right)^4.
\end{equation}

\emph{Perform in (\ref{eqxx10}) the substitution (\ref{eqc4}) which
takes the form}
\begin{equation}
y'=\frac{K(y)}{x^2}.
\end{equation}

\emph{By (\ref{eqc6}), the function} $K$  \emph{satisfies the
Bernoulli equation}
\begin{equation}\label{eqxx11}
K'_y=K+K^3
\end{equation}
\emph{ yielding the solution}
\begin{equation}
K(y)=\pm \left(  A\,e^{-2\,y}-1\right)^{-\frac{1}{2}},
\end{equation}
\emph{where} $A$ \emph{is an arbitrary constant of integration.}

\emph{Hence, by using the integral (\ref{eqc7}), the
equation (\ref{eqxx10}) provides an implicit  solution}
\begin{equation}
\arctan\left( \sqrt{A\,e^{-2\,y}-1}
\right)-\sqrt{A\,e^{-2\,y}-1}=\frac{1}{x}+B ,
\end{equation}
\emph{where} $B$ \emph{is an integration constant.}
\end{example}
\begin{example}
\emph{If} $F$ \emph{in (\ref{eqc7}) has  the
  form}
  \begin{equation}
  F(u,\,v)=v^{m+1}(h_1(u)+h_2(u)\,v+h_3(u)\,v^2),
  \end{equation}
\emph{then the corresponding reduced equation (\ref{eqc6}) leads
  to a  Riccati equation. An important number of integrable Riccati equations is recorded in [\ref{r3}-\ref{r6}]. }
\emph{Recall that, if} $y_0=y_0(x)$ \emph{is a given particular
solution of the
  Riccati equation}
  \begin{equation}
  y'=f(x)y^2+g(x)y+h(x),
  \end{equation}
  \emph{then, the general solution can be written as:}
  \begin{equation}
   y(x)=y_0(x)+\Phi(x)\left[ C-\int f(x)\Phi(x) dx
   \right]^{-1},
  \end{equation}
{\emph{where}}
  \begin{equation}
  \Phi(x)=\exp\left\{ \int\left[2f(x)y_0(x)+g(x)\right]dx  \right\};
  \end{equation}
  $C$ \emph{is an arbitrary constant. To the particular solution $y_0(x)$
  there corresponds} $C=\infty.$
\end{example}
\begin{example}
\emph{Let} $F$ \emph{in (\ref{eqc7}) be of  the
  form}
  \begin{equation}
  F(u,\,v)=v^{m+1}\frac{h_1(u)+h_2(u)\,v+h_3(u)\,v^2+h_4(u)\,v^3}{k_1(u)+k_2(u)\,v}.
  \end{equation}
  \emph{Then the corresponding reduced equation (\ref{eqc6}) leads
  to an Abel equation.}

\emph{As a matter of fact, let's briefly  present how to transform,
into its canonical form, the  Abel equations of the second kind:}
 \begin{equation}\label{eqx51}
 y'=\frac{f_2\,y^2+f_1\,y+f_0}{g_1\,y+g_0},\quad \mbox{\emph{with}}\quad f_2\neq
 0,
 \end{equation}
\emph{where} $f_i, \,i=0,\,1,\,2,$ \emph{and} $g_j,\,j=0,\,1,$
\emph{are arbitrary functions of} $x$.

\emph{Using the variable change}
\begin{equation}
y=\frac{z-g_0}{g_1},
\end{equation}
\emph{equation (\ref{eqx51}) takes the form}
\begin{equation}
z\,z'=\widetilde{f}_2\,z^2+\widetilde{f}_1\,z+\widetilde{f}_0.
\end{equation}
\emph{Now, the substitution}
\begin{equation}
z=E(x)\,w\quad \mbox{\emph{where}}\quad E(x)=e^{\int
\widetilde{f}_2(x) dx}
\end{equation}
\emph{brings this equation to the simpler form :}
\begin{equation}\label{eq14}
w\,w'_x=\widetilde{h}_1\,w+\widetilde{h}_0,
\end{equation}
\emph{which, in turn,  can be reduced, by the introduction of the
new independent variable}
\begin{equation}
s=\int \widetilde{h}_1(x) dx,
\end{equation}
\emph{into the canonical form of the second kind Abel equation }
\begin{equation}\label{eqx52}
w(s)\,w'_s=w(s)+k(s).
\end{equation}
\emph{A good compilation of integrable  Abel equations of the form
(\ref{eqx52}) can be found in [\ref{r3}-\ref{r6}].}

\emph{  Considering} $m=0,$ $a(x)=-\frac{1}{x}$ \emph{and}
$F(u,\,v)= v+2\,u, $ \emph{the equation (\ref{eqc1}) can be reduced
to the linear equation}
\begin{equation}\label{eqxxx10}
y''-\left( \frac{1}{x}+x \right)\,y'-2\,x^2\,y=0.
\end{equation}
\emph{Perform in (\ref{eqxxx10}) the substitution (\ref{eqc4}) which
takes now the form}
\begin{equation}
y'=x\,K(y).
\end{equation}
\emph{By (\ref{eqc6}), the function} $K$  \emph{verifies the Abel
equation of second kind in its canonical form}
\begin{equation}\label{eqxxx11}
K\,K'_y=K+2\,y.
\end{equation}
\emph{Two particular solutions of equation (\ref{eqxxx11}) are given
by}
\begin{equation}
K_1(y)=2\,y\quad\mbox{\emph{and}}\quad K_2(y)=-y
\end{equation}
\emph{from which can be deduced a general solution   satisfying the
algebraic equation}
\begin{equation}\label{eqxxx12}
\left(K(y)-2\,y  \right)^2\,\left( K(y)+y \right)=A,
\end{equation}
\emph{where} $A$ \emph{is an arbitrary constant.} \emph{A real
solution of the latter can be  computed to yield}
\begin{eqnarray*}
K(y)&=&\frac{1}{2}\,\sqrt[3]{-8\,y^3+4\,A+4\,\sqrt{-4\,A\,y^3+A^2}}\\
    & &
    +\frac{2\,y^2}{\sqrt[3]{-8\,y^3+4\,A+4\,\sqrt{-4\,A\,y^3+A^2}}}+y\\
    &\equiv&\Phi(y,\,A).
\end{eqnarray*}

\emph{Hence, by using the integral (\ref{eqc7}), we obtain  an implicit  solution
of the
equation (\ref{eqxxx10}) as}
\begin{equation}
\int\left[ \Phi(y,\,A) \right]^{-1}=\frac{1}{2}\,x^2+B ,
\end{equation}
\emph{where} $B$ \emph{is a constant of integration.}
\end{example}


\section{Fourth class of reducible second order  NLDEs}
Let us now investigate  second order nonlinear differential
equations of the following form:
\begin{equation}\label{eqd1}
\left(y' \right)^{m}y''+a(y)\left(y'\right)^{m+2}=e^{-(m+1)\int
a(y)dy}F\left( x,\,y'e^{\int a(y)dy}\right),
\end{equation}

\begin{remark}
\emph{It is immediate  to note that, if we set in (\ref{eqb3})}
$G(y)=C,$ \emph{where} $C$ \emph{is a constant, then equation
(\ref{eqb3}) appears as a particular case of equation (\ref{eqd1})
with} $m=0.$
\end{remark}
We first solve, as in the previous sections, the left hand side
part of equation (\ref{eqd1}), namely
\begin{equation}\label{eqd2}
\left(y' \right)^{m}y''+a(y)\left(y'\right)^{m+2}=0
\end{equation}
from which follows
\begin{equation}\label{eqd3}
\left(y' \right)^{m+1}=K^{m+1}e^{-(m+1)\int a(y) dy},
\end{equation}
where $K$ is an arbitrary constant. Suppose that $K$ is a
differentiable function of the variable $x.$ Then, equation
(\ref{eqd3}) becomes
\begin{equation}\label{eqd4}
\left(y' \right)^{m+1}=K^{m+1}(x)e^{-(m+1)\int a(y) dy}.
\end{equation}

Differentiate (\ref{eqd4}) to obtain
\begin{equation}\label{eqd5}
\left(y' \right)^{m}y''=K^{m}(x)K'_xe^{-(m+1)\int a(y) dy}
-a(y)\left(y' \right)^{m+2}.
\end{equation}

Substituting (\ref{eqd4}) and (\ref{eqd5}) into (\ref{eqd1}) we find
\begin{equation}\label{eqd6}
K^m(x)K'_x=F\left( x,\,K(x)\right).
\end{equation}

We therefore arrive at the following result:
\begin{proposition}\label{pn4}
The second order nonlinear differential equation (\ref{eqd1}) can be
reduced to the first order  differential equation (\ref{eqd6}) using
(\ref{eqd4}). Furthermore, if $K(x)=\Phi(x,\,A),$ where $A$ is an
arbitrary constant, is the general solution of (\ref{eqd6}), then
the general solution of (\ref{eqd1}) is given by
\begin{equation}\label{eqd7}
 \int e^{\int a(y) dy}dy=\int \Phi(x,\,A)dx+B,
\end{equation}
where $B$ is an arbitrary constant.
\end{proposition}

\begin{remark}
 \mbox{}
\begin{enumerate}
  \item  \emph{The equations
(\ref{eqc1}) and (\ref{eqd1}) have been investigated by Jovan  in
[\ref{r1}] \emph{for} $m=0$. Furthermore, this author has also
examined some examples of integrable equations corresponding to
particular cases where} $F$ \emph{is of the form (\ref{eqc7}) with}
$q=1,\,k_1(u)=G_1(v)=1$ \emph{and} $F_\nu(v)=v^{\alpha_\nu},$
\emph{with} $\alpha_\nu\in \IR,$ $\nu=1,\,\ldots,\,p.$

 \item \emph{The equations (\ref{eqc1}) and (\ref{eqd1}), considered  with}
 $m\neq 0,$ \emph{are equivalent to the equations of the same type with} $m=0.$
 \end{enumerate}
\end{remark}

The examples given in the previous section can be also considered in
the framework of  the Proposition \ref{pn4}. Namely, one can easily
modify those examples  to obtain equations of the form (\ref{eqd1})
which can be integrated by quadrature in simpler situations or
reduced at least to first order equations in more cumbersome cases.

As a matter of illustration, let us consider the following examples
in which the function $F$ is of the form (\ref{eqc7}) and $m=0.$

\begin{example}
\emph{Let} $a(y)=-1/y$ \emph{and}
$F(u,\,v)=\sqrt{(1-v^2)\,(1-k^2v^2)}$ \emph{where} $k\in\IR^\star.$
\emph{Equation (\ref{eqd1}) becomes}
\begin{equation}\label{eqx1}
y''-\frac{1}{y}\,(y')^2=y\,\left[ \left( 1-\left(
\frac{y'}{y}\right)^2 \right)\,\left(1-k^2\,\left( \frac{y'}{y}
\right)^2 \right) \right]^{\frac{1}{2}}.
\end{equation}

\emph{Perform in (\ref{eqx1}) the substitution (\ref{eqd4}) which
takes the form}
\begin{equation}
y'=K(x)\,y.
\end{equation}
 \emph{By (\ref{eqd6}), the function} $K$  \emph{satisfies the elliptic
integral}
\begin{equation}\label{eqx2}
\int\frac{dK}{\sqrt{(1-K^2)\,(1-k^2\,K^2)}}=x+A,
\end{equation}
\emph{where} $A$ \emph{is an arbitrary constant of integration.}
\emph{Equation (\ref{eqx2}) is readily solved to give}
\begin{equation}
K(x)=\mbox{\emph{sn}}(x+A,\,k),
\end{equation}
\emph{where sn is the first Jacobi elliptic function.} \emph{Hence,
by using the integral (\ref{eqd7}), the solution of the equation
(\ref{eqx1}) is obtained in the form }
\begin{equation}
y(x)=B\,\left[
\mbox{\emph{dn}}(x+A,\,k)+k\,\mbox{\emph{cn}}(x+A,\,k)
\right]^{-\frac{1}{k}},
\end{equation}
\emph{where} $B$ \emph{is an integration constant; cn and dn
stand for the second and  third Jacobi elliptic functions,
respectively.}
\end{example}
\section*{}
\begin{example}
\emph{Let} $a(y)=-\frac{1}{y}$ \emph{and}
$F(u,\,v)=\frac{b}{u^2}+a\,v^2,\,\,a,\,b\in \IR^\star.$
\emph{Equation (\ref{eqd1}) becomes}
\begin{equation}\label{eqxx1}
y''-\frac{1}{y}\,(y')^2=y\,\left[ \frac{b}{x^2}+a\,\left(
\frac{y'}{y} \right)^2 \right].
\end{equation}

\emph{Perform in (\ref{eqxx1}) the substitution (\ref{eqd4}) which
takes the form}
\begin{equation}
y'=K(x)\,y.
\end{equation}

\emph{By (\ref{eqd6}), the function} $K$  \emph{verifies the Riccati
equation}
\begin{equation}\label{eqxx2}
K'_x=\frac{b}{x^2}+a\,K^2
\end{equation}
\emph{which can be solved to yield}
\begin{equation}
K(x)=\frac{\lambda}{x}-x^{2\,a\,\lambda}\,\left(
\frac{a\,x^{2\,a\,\lambda+1}}{2\,a\,\lambda+1}+A \right)^{-1},
\end{equation}
\emph{where} $A$ \emph{is an arbitrary constant of integration and}
$\lambda$ \emph{is a solution of the quadratic equation}
$a\,\lambda^2+\lambda+b=0.$ \emph{Hence, by using the integral
(\ref{eqd7}), we obtain, for equation (\ref{eqxx1}), the solution}
\begin{equation}
y(x)=\frac{B\,x^\lambda}{\left( a\,x\,e^{2\,a\,\lambda\,\ln
x}+2\,a\,A\,\lambda+A \right)^{\frac{1}{a}}},
\end{equation}
\emph{where} $B$ \emph{is a constant of integration.}
\end{example}


\section{Concluding remarks}

We have investigated, in this paper,   four different  classes of
second order nonlinear  differential equations which have been
reduced to first order ones, using suitable parameter variations.
 Fortunately,  the resulting first order differential
equations are, in most cases, transformable to well known integrable or solvable classical
differential equations whose the solutions can be worked  out
 by various methods disseminated in the standard specialized text books.
Finally, it appears possible to extend the
 parameter variation methods developed in this work  to classes of higher order
  nonlinear differential equations with a view to  their order reduction.  Such an
  investigation
  will be in the core of forthcoming work.

\begin{ack}
This work is partially supported by the Abdus Salam International
Centre for Theoretical Physics (ICTP, Trieste, Italy) through the
Office of External Activities (OEA) - \mbox{Prj-15}. The ICMPA
 is in partnership with
the Daniel Iagolnitzer Foundation (DIF), France.
\end{ack}


\begin{thebibliography}{00}


\bibitem{} \label{r23} L. D. Faddeev and L. A. Takhtanjan,  \emph{ Hamiltonian methods in the theory of soliton,}
Translated from Russian by A. G. Reyman [A. G. Reiman], Springer
Series in Soviet Mathematcs. Springer-Verlag, Berlin, 1987.



\bibitem{} \label{r19} M. N. Hounkonnou and M. M. Kabir, \emph{ Hasegawa - Mima - Charney - Obukhov Equation:
Symmetry Reductions and Solutions}, Int. J. Contemp. Math. Sciences,
Vol. {\bf 3}, No. 3, (2008) 145-157.

\bibitem{} \label{r22} M. N. Hounkonnou and M. M. Kabir, \emph{ Some exact solutions
of a non linear Bousinesq system of equations}, International
Journal of pure and Applied Mathematics Vol. {\bf 45}, No. 1, (2008)
45-65.


\bibitem{} \label{r24} H. D. Doebner and G. A. Goldin, \emph{ Properties
of nonlinear Schrodinger equations associated with differomorphism
group representations}, J. Phys.  A.: Math. Gen. {\bf 27}, (1994)
1771-1780.
\bibitem{} \label{hk} M. N. Hounkonnou and M. M. Kabir, \emph{Symmetry, integrability and solutions of the
Kawahara equation}, SUT Journal of Mathematics, Vol. {\bf 44}, No.
1, (2008) 3953.
\bibitem{} \label{r20} P. J. Olver, \emph{ Application of Lie Groups to Differential
Equations},
Springer-Verlag, New York, 1993.


\bibitem{} \label{r21} L. V. Ovsiannikov, \emph{ Group Analysis of Differential
Equations},
Academic Press, New york, 1982.




\bibitem{} \label{r9} W. Malfliet, \emph{ The Tanh method, I.  Exact solutions of nonlinear
evolution and wave equations}, Physica Sprica. {\bf 54}, (1996)
569-575.


\bibitem{} \label{r10} A. M. Wazwaz,  \emph{ The Tanh method for travelling wave solutions of nonlinear
 equations}, Applied Mathematics and Computation {\bf 154}, (2004)
 713-723.

\bibitem{} \label{r11} E. Yusufoglu and A. Bekir, \emph{ Solutions of coupled nonlinear evolution
equations}, Chaos, Solitons and Fractals {\bf 37}(3), (2008)
842-848.

\bibitem{} \label{r12} E. Fan, \emph{Extended tanh-function method and its applications to nonlinear
equations}, Phys. Lett. A. {\bf 277}, (2000) 212.

\bibitem{} \label{r13} S.A. El-Wakil and M. A. Abdou, \emph{New exact travelling
wave solutions using modified extended tanh-function method}, Chaos,
Solitons and Fractals {\bf 31}(4), (2007) 840-852.

\bibitem{} \label{r14} E. Yusufoglu and A. Bekir, \emph{On the extended tanh method
applications of nonlinear equations}, International Journal of
Nonlinear Science {\bf 4}(1) (2007) 10-16.

\bibitem{} \label{r15} A. M. Wazwaz, \emph{The sine-cosine method for handling nonlinear wave
equations}, Math. and Comput. Modelling. {\bf 40}, (2004) 499-508.

\bibitem{} \label{r16} A. M. Wazwaz, \emph{The sine-cosine method for
obtaining solutions with compact and noncompact structures}, Applied
Mathematics and Computation {\bf 152}(2), (2004) 559-576.

\bibitem{} \label{r17} E. Fan and Y. C. Hon, \emph{A series of travelling wave
solutions for two variant Boussinesq equation in shallow water},
Chaos, Solitons and Fractals {\bf 15}(3), (2003) 559-566.

\bibitem{} \label{r18} Zhenya Yan, \emph{ Abundant families of Jacobi elliptic
 solutions of the (2+1)-dimensional integrable Davey-Stewartson-type equation via a new
 method},
Chaos, Solitons and Fractals {\bf 18}(2), (2003) 299-309.





\bibitem{} \label{r25} M. Wadati, \emph{ Introduction to solitons,
Pramana}, J. Phys. {\bf 57} (5-6), (2001) 841-847.

\bibitem{} \label{r26} D. LU, B. Hong and L. Tian, \emph{ Backlund transformation and
n-soliton-like solutions to the combined KdV-Burgers equation with
variable coefficients}, International Journal of Nonlinear Science
{\bf 2}, (2006) 3-10.

\bibitem{} \label{r27} M. J. Ablowitz and P. A. Clarkson, \emph{ Solitons, nonlinear
evolution equations and inverse scattering transform}, Cambridge
University Press, Cammbridge. (1990).

\bibitem{} \label{r28} P. Rosenau and J. M. Hyman, \emph{ Compactons: solitons with finite
wavelengths}, Phys. Rev. Lett., {\bf 70} (5), (1993) 564-567.

\bibitem{} \label{r29} G. Cai, Q. Wang and J. Huang, \emph{ A modified
 f-expansion method for solving breaking soliton equation},
International Journal of Nonlinear Science {\bf 2}, (2006) 122-128.

\bibitem{} \label{r30} D. Zhang, \emph{ Doubly periodic solutions of modified Kawahara
equation}, Chaos, Solitons and Fractals {\bf 26}, (2005) 1155-1160.

\bibitem{} \label{r31} H. Zhang, \emph{ New exact travelling wave solutions for some nonlinear
evolution equations}, Chaos, Soliton and Fractals {\bf 25}, (2005)
921-925.

\bibitem{} \label{r1} Jovan D. Ke\v{c}ki\'{c}, \emph{Additions to Kamke's treatise,
 VII : Variation  of parameters
for nonlinear second order differential equations}, Univ. Beograd.
Pool. Elektrotehn. fak. Ser. Mat. Fiz. No. {\bf 544} - No. {\bf
576}, (1946) 31-36.




\bibitem{} \label{r2} E. Kamke, \emph{Differentialgleichungen: Losungsmethoden und Losungen, I,
Gewohnliche Differentialgleichungen}, B. G. Teubner, Leipzig, 1977.



\bibitem{} \label{r3} G. M. Murphi, \emph{Ordinary Differential Equations and Their Solutions},
D. Van Nostrand, New York, 1960.



\bibitem{} \label{r4} A. D. Polyanin and V. F. Zaitsev, \emph{Handbook of Exact Solutions for Ordinary Differential Equations},
Chapman and Hall/CRC Press, Boca Raton, 2nd edition 2003.



\bibitem{} \label{r5} D. Zwillinger, \emph{Handbook of Differential Equations},
Academic Press, Boston, 3rd edition, 1997.




\bibitem{} \label{r6} A. D. Polyanin and A. V. Manzhirov, \emph{Handbook of Mathematics for Engineers and Scientists}
(Chapters 12, T5, and T6), Chapman and Hall/CRC Press, Boca Raton
London, 2006.



\bibitem{} \label{r7} W. E. Boyce and R. C. DiPrima,  \emph{Elementary Differential Equations,
7th Edition}, Wiley, New York, 2000.



\bibitem{} \label{r8} P. Appell and E. Lacour: \emph{Principe de la theorie des functions elliptiques et applications.}
Paris, Gauthier - Villars et Fils, imprimeurs - libraires de l'école
polytechnique, du bureau des longitudes, 23089 Quai des Grands
Augustins, 55, 1897.











\end{thebibliography}
\end{document}